\documentclass[11pt]{amsart}

\usepackage{amssymb,amsmath,amsxtra,amscd,eucal}

\theoremstyle{plain}
\numberwithin{equation}{section}
\newtheorem{theorem}{Theorem}

\numberwithin{lemma}{section}
\numberwithin{theorem}{section}
\numberwithin{corollary}{section}
\numberwithin{proposition}{section}

\theoremstyle{definition}
\newtheorem{definition}{Definition}
\numberwithin{definition}{section}
\theoremstyle{remark}

\numberwithin{example}{section}
\newtheorem{remark}{Remark}

\newcommand{\nc}{\newcommand} 
\nc{\mc}{\mathcal}
\nc{\on}{\operatorname} 
\nc{\BunG}{\on{Bun}_{G}}
\nc{\BunGH}{\on{\BunG^{H}}} 
\nc{\Mg}{\on{M}^{g}} 
\nc{\M}{\on{M}}
 
\nc{\Cl}{Cl} 
\nc{\drva}{\widehat{\Omega}}
\nc{\spec}{\on{Spec}} 
\nc{\vac}{|0\rangle}
\nc{\Z}{\mathbb{Z}} 
\nc{\zf}[1]{z^{\frac{1}{#1}}}
\nc{\wf}[1]{w^{\frac{1}{#1}}} 
\nc{\Mt}{M^{\sigma}}
\renewcommand{\Mg}[1]{\mathfrak{M}_{g,#1}} 
\nc{\MGg}[1]{\mathfrak{M}^{G}_{g,#1}}
 \nc{\gr}{\on{gr}}

\nc{\fps}[1]{\ensuremath{\mathbb{C}[[#1]]}}
\nc{\ls}[1]{\ensuremath{\mathbb{C}((#1))}}
\nc{\tf}[1]{\ensuremath{t^{\frac{1}{#1}}}}
\newcommand{\ru}[1]{\ensuremath{e^{\frac{2 \pi i}{N}}}}

\nc{\Orb}[1]{\ensuremath{\mathbf{O}_{#1}}}

\nc{\Cnot}{\ensuremath{C \backslash \{ \nu^{-1}(x_{i}) \}}}

\renewcommand{\k}{\ensuremath{\mathfrak{k}}}
\nc{\g}{\ensuremath{\mathfrak{g}}}
\renewcommand{\l}{\ensuremath{\mathfrak{l}}}
\nc{\goutc}{\ensuremath{\g^{\sigma}_{\on{out}}(\Caff)}}
\nc{\Lg}{\ensuremath{L \g}}
\nc{\ghat}{\ensuremath{\hat{\g}}}

\nc{\Ksto}{\ensuremath{\mathcal{K}^{*}_{odd}}}
\nc{\Osto}{\ensuremath{\mathcal{O}^{*}_{odd}}}
\nc{\Ost}{\ensuremath{\mathcal{O}^{*}}}

\nc{\Vkg}{\ensuremath{V_{k}(\mathfrak{g})}}

\nc{\V}{\ensuremath{\mathcal{V}}}
\nc{\Vt}{\ensuremath{\wt{\mathcal{V}}}}

\nc{\nV}{\nu^*(\V)}

\nc{\Vtx}{\ensuremath{\widetilde{\mathcal{V}}_{X}}}
\nc{\Ptx}{\ensuremath{\Pi_{C}}}
\nc{\Mtp}[1]{\ensuremath{\mathcal{M}^{\sigma}_{#1}}}
\nc{\Mp}[1]{\ensuremath{\mathcal{M}_{#1}}}
\nc{\Np}[1]{\ensuremath{\mathcal{N}_{#1}}}
\nc{\Mtc}{\ensuremath{\mathcal{M}}}
\nc{\Ym}[1]{\ensuremath{Y^{#1}}}

\nc{\fd}[1]{\ensuremath{\mathcal{D}_{#1}}}
\nc{\fpd}[1]{\ensuremath{\mathcal{D}^{\times}_{#1}}}
\nc{\AutO}{\ensuremath{\Aut ( \mathcal{O} )} }
\nc{\AutOr}[1]{\Aut_{#1} \mathcal{O}}

\nc{\AutD}[1]{\Au(\mathcal{D}_{#1})}
\nc{\AutDN}[2]{\Au_{#1}(\mathcal{D}_{#2})}
\nc{\AutX}[1]{\ensuremath{\emph{\mc{A}ut}_{#1}}}

\nc{\discs}[1]{\ensuremath{\coprod_{p \in \nu^{-1}(x_{#1})} \fpd{p}}}

\nc{\vect}{\ensuremath{\mc{V} ect}}

\nc{\Spec}{\on{Spec}}
\nc{\Aut}{\ensuremath{\on{Aut}}}
\nc{\Der}{\on{Der}}
\nc{\End}{\on{End}}
\nc{\Au}{\ensuremath{ {\mc A}ut}}
\nc{\ol}[1]{\ensuremath{\overset{\circ}{#1}}}
\nc{\ep}{\epsilon}
\nc{\C}{{\mathbb C}}
\nc{\OO}{{\mc O}}
\nc{\K}{{\mc K}}
\nc{\DD}{{\mc D}}
\nc{\bD}{\overline{\DD}}
\nc{\pa}{\partial}
\nc{\secref}[1]{Section~{\ref{#1}}}
\nc{\corref}[1]{Corollary~{\ref{#1}}}
\nc{\wt}{\widetilde}
\nc{\Res}{\on{Res}}
\nc{\Ind}{\on{Ind}}
\nc{\nn}{_{(n)}}
\nc{\cA}{{\mc A}}
\nc{\Wick}{{\mathbf :}}
\nc{\Hom}{\on{Hom}}
\nc{\al}{\alpha}
\nc{\mb}{\mathbf}
\nc{\ds}{\displaystyle}

\begin{document}

\title{Orbifold conformal blocks and the stack of pointed G-covers}

\author{Matthew Szczesny}
\address{Department of Mathematics  
         University of Pennsylvania, Philadelphia, PA 19104 }
\email{szczesny@math.upenn.edu}

\begin{abstract}
Starting with a vertex algebra $V$, a finite group $G$ of
automorphisms of $V$, and a suitable collection of twisted
$V$--modules, we construct (twisted) $D$--modules on the stack of
pointed $G$--covers, introduced by Jarvis, Kaufmann, and Kimura. The
fibers of these sheaves are spaces of orbifold conformal blocks
defined in joint work with Edward Frenkel. The key ingredient is a
$G$--equivariant version of the Virasoro uniformization theorem.
\end{abstract}

\maketitle 


\section{Introduction}

It is by now well-understood that conformal field theories (CFT's) in
two dimensions give rise to sheaves with projectively flat connection
over the moduli stack of $n$-pointed genus $g$ curves,
$\Mg{n}$. Mathematically, this process can be described as follows
(see for instance the book \cite{FB} for more details). Given a vertex
algebra $V$, which corresponds to a choice of CFT model, a collection
of $V$--modules
\[
 \mathbb{M} = (M_{1}, \cdots, M_{n}),
\]
and an n-pointed curve $(X, {\bf p}),\; {\bf p} = (p_{1}, \cdots,
p_{n})$, one obtains the vector space of conformal blocks, denoted
$C_{V}(X, \bf{p}, \mathbb{M})$. Elements of this vector space can be
used to construct chiral correlation functions in the CFT, which are
sections of certain sheaves on powers of $X$ with pairwise diagonals
removed.  As $(X, \bf{p}) $ varies in $\Mg{n}$, one obtains a sheaf
with a projectively flat connection, known as the
Knizhnik-Zamolodchikov (KZ) connection. In certain cases, when the CFT is
rational, the spaces $C_{V}(X, \bf{p}, \mathbb{M})$ are
finite-dimensional.  The sheaf is then a vector bundle that extends to
the Deligne-Mumford compactification $\overline{\Mg{n}}$, and the
connection to one with logarithmic singularities along the
boundary divisor.

In several cases, the space of conformal blocks is related to a moduli
problem. For instance, when $V = L_{k}(\g)$, the integrable basic
$\ghat$--module of level $k \in \mathbb{Z}_{+}$, we have the
isomorphism
\[
C_{L_{k}(\g)}(X,p,L_{k}(\g)) \cong H^{0}( \BunG (X), \Theta^{k})
\]
where $\BunG (X)$ is the moduli stack of $G$-bundles on $X$, $G$ is
the simply-connected algebraic group with Lie algebra $\g$, and
$\Theta$ is the theta line bundle on $\BunG (X)$, In this case, the
KZ connection yields a non-abelian analogue of the heat equation
satisfied by abelian theta functions.

In recent years, much attention has been paid to the construction of
CFT's by orbifolding. At the level of vertex algebras, this procedure
can be described as follows.  We are given a vertex algebra $V$, and a finite
group $G$ of automorphisms of $V$ that also acts on the category of
$V$--modules. One begins by taking $G$--invariants, and then adjoining
so called twisted modules $\{ M^{[g]} \}$, for each conjugacy class $[g]$
in $G$. Geometrically, orbifold models possess the new feature of
twist fields. These are objects that cause vertex operators to become
multivalued, with monodromies given by elements of $G$. 

In \cite{FS}, the notion of conformal block was extended to include
twisted modules. A key ingredient is the notion of a pointed
$G$--cover. Let $(X, {\bf p})$ be a an $n$--pointed
smooth projective curve. A pointed $G$--cover of $(X, {\bf p})$ is a set of data 
$$ 
(\pi:C \mapsto X, {\bf p}, {\bf q}), \; \; {\bf q} = (q_{1}, \cdots, q_{n}) 
$$
where $C$ is a smooth projective curve (not
necessarily connected) with an effective action of $G$, such that $X =
C/G$, the quotient map $\pi: C \mapsto X$ makes $C$ in to a principal
$G$--bundle over $X \backslash \{ p_{i} \}$, and $q_{i} \in
\pi^{-1}(p_{i})$. The space of orbifold conformal blocks, denoted
\begin{equation} \label{blocks}
C^{G}_{V}(C,X, {\bf p}, {\bf q}, \mathbb{M})
\end{equation}
is attached to a $G$--cover, and a collection of $V$--modules
$\mathbb{M} = (M_{1}, \cdots, M_{n})$, where $M_{i}$ is twisted by the
monodromy generator $m_{i}$ at $q_{i}$. Orbifold conformal blocks can
be used to construct chiral orbifold correlators, which are
$G$--invariant sections of sheaves on powers of $C$ with certain
divisors removed. For a detailed treatment of this construction, see
our forthcoming paper \cite{S}.

In \cite{JKK}, following earlier work by \cite{ACV}, the authors
introduce a smooth Deligne-Mumford stack $\MGg{n}$ parameterizing
smooth $n$--pointed $G$--covers. Given an appropriate collection
$\mathbb{M}$ of (twisted) $V$--modules, it is a natural question how
the spaces \ref{blocks} vary as $(\pi:C \mapsto X, {\bf p}, {\bf q})$
moves in $\MGg{n}$.  In this paper, we construct sheaves with a projectively
flat connection over (certain components of) $\MGg{n}$, whose fibers
at a $G$--cover $(\pi:C \mapsto X, {\bf p}, {\bf q}) $ are the
corresponding spaces \ref{blocks}.  In other words, we construct a
localization functor
\[
\Delta: \left( V- \on{mod} \right) \longrightarrow  \left( \wt{D}_{\MGg{n}}- \on{mod} \right)
\] 
where $ \wt{D}_{\MGg{n}}$ is a sheaf of twisted differential
operators on $\MGg{n}$, depending on the Virasoro central charge. This
yields an orbifold generalization of the KZ connection. The resulting
sheaves should be useful in constructing G-modular functors (see
\cite{Ki}).

In the untwisted case, the key theorem in the construction is the so
called Virasoro uniformization of $\Mg{n}$ (see \cite{ADKP},
\cite{BS}, \cite{Kon}, \cite{TUY} ). Roughly, this theorem states that
$\Mg{n}$ is a flag manifold for the Virasoro algebra. We extend the
approach to $G$--covers.

Just as $\Mg{n}$ can be compactified to $\overline{\Mg{n}}$, so does
$\MGg{n}$ have a compactification $\overline{\MGg{n}}$. In this paper
we do not address the behavior of the $D$--modules on the
boundary $\overline{\MGg{n}} \backslash \MGg{n}$. Just as in the
untwisted case, the existence of a logarithmic extension requires the finite-dimensionality of conformal
blocks. It is natural to conjecture that such an extension is possible
for the WZW model of positive integer level, when $V= L_{k}(\g)$, and
where the twisted modules are realized as integrable representations
of appropriate twisted affine algebras.

The structure of the paper is as follows. In section \ref{VOA} we
quickly review the notion of vertex algebra and vertex algebra
module/twisted module. Section \ref{Gcovers} recalls some facts about
the stack of pointed $G$--covers from \cite{JKK}. In sections
\ref{Torsors} and \ref{ConfBlocks} we review the construction of the
space of orbifold conformal blocks $C^{G}_{V}(C,X, {\bf p}, {\bf q},
\mathbb{M})$ following \cite{FS}.  Section \ref{Localization} is
devoted to reviewing Beilinson-Bernstein localization, which is the
technique used to produce $D$--modules on $\MGg{n}$. Our
treatment is essentially a condensed version of that in
\cite{FB}. Section \ref{Uniformization} contains a proof of the
$G$--equivariant Virasoro uniformization of $\MGg{n}$ and the
construction of the localization functor. Finally, section \ref{FixedBase} gives a formula for the orbifold KZ connection along the fibers of the projection 
\[
\MGg{n} \mapsto \Mg{}
\]
sending $(\pi: C \mapsto X, {\bf p}, {\bf q})$ to $X$, i.e. the locus of pointed $G$--covers with a fixed base curve. 

\bigskip

\noindent {\bf Acknowledgments:} I would like to thank Edward Frenkel
and David Ben-Zvi for many valuable conversations. I would also like to greatfully acknowledge the hospitality of  the Erwin Schrodinger Institute, where part of this work was completed. This research was
partially supported by NSF grant DMS-0401619.

\section{Vertex algebras and twisted modules} \label{VOA}

In this paper we will use the language of vertex algebras, their
modules, and twisted modules. For an introduction to vertex algebras
and their modules see \cite{FLM,K,FB}, and for background on twisted
modules, see \cite{FFR,Dong,DLM, BK}.

We recall that a conformal vertex algebra is a $\Z_+$--graded vector
space $$V = \bigoplus_{n=0}^\infty V_n,$$ together with a vacuum
vector $\vac \in V_0$, a translation operator $T$ of degree $1$, a
conformal vector $\omega \in V_2$ and a vertex operation
\begin{align*}
Y: V &\to \on{End} V[[z^{\pm 1}]], \\
A &\mapsto Y(A,z) = \sum_{n \in \Z} A_{(n)} z^{-n-1}.
\end{align*}
These data must satisfy certain axioms (see \cite{FLM,K,FB}). In
what follows we will denote the collection of such data simply by $V$.

A vector space $M$ is called a $V$--module if it is equipped with an
operation
\begin{align*}
Y^M: V &\to \on{End} M[[z^{\pm 1}]], \\
A &\mapsto Y^M(A,z) = \sum_{n \in \Z} A^M_{(n)} z^{-n-1}
\end{align*}
such that for any $v \in M$ we have $A^M_{(n)} v = 0$ for large enough
$n$. This operation must satisfy the following axioms:

\begin{itemize}

\item $Y^M(\vac,z) = \on{Id}_M$;

\item For any $v \in M$ there exists an element $$f_v \in
M[[z,w]][z^{-1},w^{-1},(z-w)^{-1}]$$ such that the formal power series
$$Y^M(A,z) Y^M(B,w) v \qquad \on{and} \qquad Y_M(Y(A,z-w) B,w) v$$ are
expansions of $f_v$ in $M((z))((w))$ and $M((w))((z-w))$,
respectively.

\end{itemize}

The power series $Y^M(A,z)$ are called vertex operators. We write the
vertex operator corresponding to $\omega$ as
\[
	Y^M(\omega,z) = \sum_{n \in \mathbb{Z}} L^M_{n} z^{-n-2},
\]
where $L^M_n$ are linear operators on $V$ generating the Virasoro
algebra. Following \cite{Dong}, we call $M$ \emph{admissible} if
$L^{M}_{0}$ acts semi-simply with integral eigenvalues.

Now let $\sigma_{V}$ be a conformal automorphism of $V$, i.e., an
automorphism of the underlying vector space preserving all of the
above structures (so in particular $\sigma_{V}(\omega) = \omega$). We
will assume that $\sigma_{V}$ has finite order $N>1$. A vector space
$M^\sigma$ is called a $\sigma_V$--{\em twisted} $V$--module (or
simply twisted module) if it is equipped with an operation
\begin{align*}
Y^{M^\sigma}: V &\to \on{End} M^\sigma[[z^{\pm \frac{1}{N}}]], \\ A
&\mapsto Y^{M^{\sigma}}(A,\zf{N}) = \sum_{n \in \frac{1}{N}\Z}
A^{M^\sigma}_{(n)} z^{-n-1}
\end{align*}
such that for any $v \in M^{\sigma}$ we have $A^{M^\sigma}_{(n)} v =
0$ for large enough $n$. Please note that we use the notation
$Y^{M^{\sigma}}(A,\zf{N})$ rather than $Y^{M^{\sigma}}(A,z)$ in the
twisted setting. This operation must satisfy the following axioms (see
\cite{FFR,Dong,DLM,LI}):

\begin{itemize}

\item $Y^{M^{\sigma}}(\vac,\zf{N}) = \on{Id}_{\Mt}$;

\item For any $v \in \Mt$, there exists an element 
\[
f_{v} \in \Mt [[\zf{N}, \wf{N} ]][z^{- \frac{1}{N}}, w^{-
\frac{1}{N}},(z-w)^{-1}]
\]
such that the formal power series $$Y^{\Mt}(A,\zf{N}) Y^{\Mt}(B,
\wf{N})v \qquad \on{and} \qquad Y^{\Mt}(Y(A,z-w) B, \wf{N})v$$ are
expansions of $f_{v}$ in $\Mt((\zf{N}))((\wf{N}))$ and
$\Mt((\wf{N}))((z-w))$, respectively.

\item If $A \in V$ is such that $\sigma_V(A) = e^{\frac{2\pi i m}{N}}
A$, then $A^{M^\sigma}_{(n)} = 0$ unless $n \in \frac{m}{N} + \Z$.

\end{itemize}

The series $Y^{M^\sigma}(A,z)$ are called twisted vertex operators.
In particular, the Fourier coefficients of the twisted vertex operator
\[
	Y^{M^\sigma}(\omega,\zf{N}) = \sum_{n \in \mathbb{Z}}
	L^{M^\sigma}_{n} z^{-n-2},
\]
generate an action of the Virasoro algebra on $M^\sigma$. The
$\sigma_{V}$--twisted module $\Mt$ is called \emph{admissible} if
$L^{\Mt}_{0}$ acts semi-simply with eigenvalues in $\frac{1}{N}\Z$.

Suppose that $\Mt$ is an admissible module. Then we define a linear
operator $S_\sigma$ on $\Mt$ as follows. It acts on the eigenvectors
of $L^{\Mt}_0$ with eigenvalue $\frac{m}{N}$ by multiplication by
$e^{\frac{2\pi i m}{N}}$. Hence we obtain an action of the cyclic
group of order $N$ generated by $\sigma$ on $\Mt$, $\sigma \mapsto S_\sigma$.
According to the axioms of twisted module,  we have the following identity:
\begin{equation}    \label{action on module}
S_\sigma^{-1} Y^{\Mt}(\sigma \cdot A,\zf{N}) S_\sigma = Y^{\Mt}(A,\zf{N}).
\end{equation}

\section{The stack of pointed G-covers} \label{Gcovers}

In this section we review the definition of the stack of pointed
G-covers introduced in \cite{JKK} following \cite{ACV}.

Let $\Mg{n}$ denote the stack of smooth $n$-pointed curves of genus
$g$. The objects of $\Mg{n}$ consist of families $(\lambda: X \mapsto
S, p_{1}, \cdots, p_{n})$ where $\lambda : X \mapsto S$ is a smooth
family of curves of genus $g$, and $p_{i}: S \mapsto X, \; i=1,
\cdots, n$ are pairwise disjoint sections of $\lambda$. Note that the
$p_{i}$ are ordered.

\begin{definition}
Let $(\lambda: X \mapsto S, p_{1}, \cdots, p_{n}) \in \Mg{n}$ be a
smooth $n$--pointed curve, and $G$ a finite group.  A smooth
$G$--cover of $X$ consists of a morphism $\pi: C \mapsto
X$ of smooth curves over $S$, satisfying the following properties:
\begin{itemize}
	\item There is a right $G$--action on $C$ preserving $\pi$. 
	\item $C \backslash \cup_{i} \pi^{-1}(p_{i})$ is a principal
          $G$--bundle over $X \backslash \cup_{i} p_{i}$
\end{itemize}
We denote the smooth $G$--cover by $(\pi: C \mapsto X, p_{1}, \cdots ,
p_{n})$. A smooth $n$--pointed $G$--cover consists of a smooth
$G$--cover $(\pi:C \mapsto S, p_{1}, \cdots, p_{n})$ together with $n$
sections $q_{i}: S \mapsto C$ such that $q_{i} \in
\pi^{-1}(p_{i})$. We denote it by $(\pi: C \mapsto X, p_{1}, \cdots,
p_{n}, q_{1}, \cdots, q_{n})$.
\end{definition} 

To avoid cumbersome notation, we will henceforth use $(X,
{\bf p})$, $(C,X,{\bf p})$, and $(C,X,{\bf p}, {\bf q})$, to denote
respectively, an $n$--pointed curve, a $G$--cover, and a pointed
$G$--cover, where ${ \bf p } = (p_{1}, \cdots, p_{n})$ and ${\bf q} =
(q_{1}, \cdots, q_{n})$.

A morphism of pointed $G$--covers is a $G$--equivariant fibered
diagram - that is, a morphism of the underlying curves $X$ together
with a $G$--equivariant morphism of the covers preserving the points
$q_{i}$.  Smooth $n$--pointed $G$--covers form a stack, denoted
$\MGg{n}$. There is an obvious morphism
\begin{align*}
p: \MGg{n} & \mapsto \Mg{n} \\ (C,X,{\bf p}, {\bf q}) & \mapsto (X, {\bf p})
\end{align*}

Let $G_{A}$ denote the group $G$ viewed as a right $G$--space under
conjugation. There is a $G$--equivariant map
\[
{\bf e}:  \MGg{n} \mapsto G^{n}_{A}
\]
sending $(C,X, {\bf p},  {\bf q})$ to the $n$--tuple ${\bf m} = (m_{1}, \cdots, m_{n})$, where $m_{i}$ is
the monodromy of $C$ around $q_{i}$. Let $\MGg{n}(\bf{m})$ denote the
closed sub-stack ${\bf e}^{-1}(\bf{m})$. Note that $\MGg{n}(\bf{m})$
may be empty. We have
\[ 
\MGg{n} = \bigcup_{{\bf m} \in G^{n}_{A} } \MGg{n}(\bf{m})
\]
We have the following theorem:

\begin{theorem}{\cite{JKK}}
The stack $\MGg{n}$ and the stacks $\MGg{n}(\bf{m})$ are smooth
Deligne-Mumford stacks, flat, proper, and quasi-finite over $\Mg{n}$.
\end{theorem} 

\section{Groups and torsors} \label{Torsors}

The purpose of this section is to review the structure of the group of
formal coordinate changes on the "formal" disk $\Spec(R[[z]])$ and the
"formal punctured disk" $\Spec(R((z)))$ for an arbitrary
$\mathbb{C}$--algebra $R$. We will be primarily interested in the case
when $R$ is an Artin $\mathbb{C}$-algebra, such as
$\mathbb{C}[\epsilon]/(\epsilon^{2})$. For more on the structure of
these groups, see Section 5.1 of \cite{FB}, or the original source
\cite{MP}.

\subsection{Groups}

Let us denote $R[[z]]$ by $\OO_{R}$, and $R((z))$ by $\mc{K}_{R}$. Let
$\Aut (\mc{O}_{R})$ denote the group of continuous $R$--algebra
automorphisms of $R[[z]]$ preserving the ideal $(z)$, and
$\Aut(\mc{K}_{R})$ the group of all continuous $R$--algebra
automorphisms of $R((z))$. Since $R[[z]]$ is topologically generated
by $z$, an automorphism $\rho$ of $R[[z]]$ is completely determined by
the image of $z$, which is a series of the form
\begin{equation}
\rho(z) = \sum_{n \in \mathbb{Z}, n \geq 1}
c_{n} z^{n},  \qquad c_{n} \in R  \label{rhodef1}
\end{equation}
where $c_{1}$ is a unit. Hence we identify $\Aut (\mc{O}_{R})$
with the space of power series in $z$  satisfying these conditions.  
Similarly, an element $\rho$ of $\Aut(\mc{K}_{R})$ can be identified
with a formal Laurent series of the form
\[
\rho(z) = \sum_{n \in \mathbb{Z}, n \geq k}
c_{n} z^{n},
\]
where $c_{n}$ is nilpotent if $n \leq 0$, and $c_{1}$ is a
unit. The functor
\[
   R \mapsto \Aut(\mc{O}_{R})
\]
is representable by a group scheme over $\mathbb{C}$ which we'll
denote $\Aut \mc{O}$.  The Lie algebra $\Der^{(o)}(\mc{O})$ of $\Aut(\mc{O}
)$ is topologically generated by elements
\[
z^{k} \partial_{z}, \; k \geq 1
\]
The functor
\[
 R \mapsto \Aut(\mc{K}_{R})
\]
is representable by an Ind-group scheme which we'll denote $\Aut(\mc{K})$. 
Denote by $\Der(\mc{K})$ the Lie algebra of $\Aut(\mc{K})$, generated by
\[
 z^{k} \partial_{z}, \; k \in \mathbb{Z} 
\]

We now consider some groups arising in the study of ramified coverings of disks.
Let $\Aut  (R[[\zf{N}]])$ denote the group of continuous $R$--algebra
automorphisms of $R[[\zf{N}]]$ preserving the ideal $(\zf{N})$.

\begin{definition} $\Aut_{N} (\mc{O}_{R})$ is the subgroup of $\Aut (R[[\zf{N}]])$
preserving the subalgebra $R[[z]] \subset R[[z^{\frac{1}{N}}]]$.
\end{definition}

\noindent Thus, $\AutOr{N}$ consists of power series of the form
\begin{equation}
\rho(\zf{N}) = \sum_{n \in \frac{1}{N} + \mathbb{Z}, n > 0}
c_{n} z^{n}, \qquad c_{n} \in R  \label{rhodef}
\end{equation}
such that $c_{\frac{1}{N}}$ is a unit. 

There is a homomorphism $\mu_{N}: \Aut_{N} (\mc{O}_{R}) \rightarrow
\Aut (\mc{O}_{R})$ which takes $\rho \in \Aut_{N} (\mc{O}_{R})$ to the
automorphism of $R[[z]]$ that it induces. At the level of power
series, this is just the map $\mu_{N}: \rho(\zf{N}) \mapsto
\rho(\zf{N})^{N}$. The kernel consists of the automorphisms of the
form $\zf{N} \mapsto \epsilon \zf{N}$, where $\epsilon$ is an $N$th
root of unity, so we have the following exact sequence:
\[
	1 \rightarrow \mathbb{Z}/ N \mathbb{Z} \rightarrow
	\Aut_{N} (\mc{O}_{R}) \rightarrow \Aut(\mc{O}_{R}) \rightarrow 1 \, .
\]
Thus $\Aut_{N} (\mc{O}_{R})$ is a central extension of
$\Aut(\mc{O}_{R})$ by the cyclic group $\mathbb{Z}/ N \mathbb{Z}$. One
can define the group $\Aut_{N}(\mc{K}_{R})$ in an obvious way. Denote
by $\Aut_{N}(\mc{O})$ the group scheme representing the functor $R
\mapsto \Aut_{N}(\mc{O}_{R})$, and $\Aut_{N}(\mc{K})$ the Ind-group scheme
representing $R \mapsto \Aut_{N}(\mc{K}_{R})$.
The Lie algebra $\Der^{(o)}_{N} (\mc{O})$ of $\Aut_{N}(\mc{O})$
 can be identified with 
 $\zf{N} \mathbb{C}[[z]] \partial_{\zf{N}}$. The homomorphism $\mu_{N}$ induces an
isomorphism of the corresponding Lie algebras sending
\begin{equation} \label{muN}
	z^{k+\frac{1}{N}} \partial_{\zf{N}} \mapsto N z^{k+1}
	\partial_{z}, \qquad k \in \mathbb{Z}, k \geq 0.
\end{equation}
Similarly, the Lie algebra $\Der_{N}(\mc{K}) $of $\Aut_{N}(\mc{K})$
can be identified with $\zf{N} \mathbb{C}((z)) \partial_{\zf{N}}$, and
$\mu_{N}$ extends to a homomorphism $\Der_{N}(\mc{K}) \mapsto \Der(\mc{K})$.


Let $X$ be a smooth curve over $\Spec R$. If $x \in X$ is an
$R$--valued point of $X$, denote by $\Aut(\widehat{\mc{O}}_{x,R})$ the
automorphisms of the formal neighborhood of $x$ fixing $x$.  Suppose
now that $C$ is a smooth curve over $\Spec R$ with effective
$G$--action, and that $X = C/G$. Denote the quotient map by $\pi: C
\mapsto X$.
Choosing a formal coordinate $z$  at
$x$ yields an isomorphism  $\Aut (\widehat{\mc{O}}_{x,R}) \cong \Aut (\mc{O}_{R})$ (For the definition of a formal coordinate, see the next section).

If $y \in C$, denote by $O(y)$ the $G$--orbit of $y$ in $C$, and by
$\widehat{O(y)}$ the formal neighborhood of $O(y)$. 
Let $\widehat{O(y)}^{*}$ denote the union of the formal
punctured disks around points in $O(y)$.
Note that $\widehat{O(y)}$ and $\widehat{O(y)}^{*}$ carry an action of $G$.  Let
\[
\Aut(\widehat{O(y)}_{R}) \; resp. \Aut(\widehat{O(y)}^{*}_{R})
\]
 denote the automorphism group of $O(y)$ and $O(y)^{*}$ respectively, and 
\[
\Aut^{G}(\widehat{O(y)}_{R}) \; resp. \Aut^{G}(C,
\widehat{O(y)}^{*}_{R})
\]
the corresponding subgroups of elements commuting with $G$. 
Finally, denote by 
\[
\Aut^{G}_{e}(\widehat{O(y)}_{R}) \; resp. \Aut^{G}_{e}(
\widehat{O(y)}^{*}_{R})
\]
the identity components of the corresponding groups.  For each $q \in
O(y)$, by choosing a coordinate in which the action of the stabilizer
$G_{q}$ is linear (this is called a special coordinate, see the next
section), we obtain isomorphisms
\[
\psi_{q}: \Aut^{G}_{e} (\widehat{O(y)}_{R}) \mapsto  \Aut_{N} (\mc{O}_{R}) 
\]
and
\[
\psi^{*}_{q}: \Aut^{G}_{e} ( \widehat{O(y)}^{*}_{R}) \mapsto  \Aut_{N} (\mc{K}_{R}) 
\]
by restricting to the disk around $q$. The proof amounts to observing
that any element of $\Aut_{N} (\mc{O}_R)$ has a unique
$G$--equivariant extension to $ \Aut^{G}_{e}(\widehat{O(y)}_{R})$ 
(and the same with $\widehat{O}$ replaced by $\widehat{O}^{*}$).

\subsection{Torsors}

Let $\mathcal{D}_{R} = \Spec A$, where $A \cong R[[z]]$ is a 
``formal'' $R$--disk, and let $x \in \mathcal{D}_{R}$ be an
$R$--point.  By a formal coordinate on $(\mc{D}_{R},x)$ we mean a
isomorphism $\mc{D}_{R} \cong \Spec(R[[z]])$ that identifies $x$ with the
origin (the origin being the $R$--point corresponding to the ideal
$(z) \subset R[[z]]$). Let $\Au(\mc{D}_{R,x})$ denote the set of
formal coordinates on $(\mc{D}_{R},x)$. It is an
$\Aut(\mc{O}_{R})$--torsor.

Suppose now that $(\mathcal{D}_{R}, x, \sigma_{\mathcal{D}})$ is a
triple consisting of a formal disk $\mathcal{D}_{R} = \Spec B$, where
$B \cong R[[\zf{N}]]$, $x \in \mc{D}_{R}$, and $\sigma_{\mathcal{D}}$
is an automorphism of $\mathcal{D}$ (equivalently, of $B$) of order
$N$ fixing $x$. After a change of coordinate, $\sigma_{\mc{D}}$ is
equivalent to the automorphism $\zf{N} \mapsto \epsilon \zf{N}$, where
$\epsilon$ is a primitive $N$th root of unity.  We denote by $\bD$ the
quotient of $\DD$ by $\langle \sigma_{\DD} \rangle$, i.e., the disk
$\Spec B^{\sigma_{\DD}}$, where $B^{\sigma_{\DD}}$ is the subalgebra
of $\sigma_{\DD}$--invariant elements.

A formal coordinate $t$ is called a {\em special coordinate} with
respect to $\sigma_{\mc{D}}$ if $\sigma_{\mc{D}}(t) = \epsilon t$,
where $\epsilon$ is an $N$th root of unity, or equivalently, if $t^N$
is a formal coordinate on $\bD$. We denote by $\Au_{N}(\mc{D}_{R,x})$
the subset of $\Au(\mc{D}_{R,x})$ consisting of special formal
coordinates. The set $\Au_{N}(\mc{D}_{R,x})$ carries a simply
transitive right action of the group $\Aut_{N}(\mc{O}_{R})$ given by
$t \mapsto \rho(t)$, where $\rho$ is the power series given in
\eqref{rhodef}, i.e. $\Au_{N}(\mc{D}_{R,x})$ is an
$\Aut_{N}(\mc{O}_{R})$--torsor.

Suppose now that $\pi: C \mapsto X=C/G$ as above, and that $y \in
C$. Let $N$ denote the order of the stabilizer $G_{y}$, which is cyclic, and
generated by the monodromy around $y$, $m_{y}$. To this data we can
associate the set $\Au_{N}(\mc{D}_{R,y})$ of special coordinates on 
$(\Spec \widehat{\mc{O}}_{y}, y)$ with respect to $G_{y}$,
which is an $\Aut_{N}(\mc{O}_{R})$--torsor. 

Note that when $R=\mathbb{C}$, $\mc{D}_{R}$ has a unique $R$--point, and we will suppress it in the notation. Henceforth, we will also use the convention that $\mc{D}_{\mathbb{C}}$ is denoted $\mc{D}$,
and suppress $R$ when referring to groups, torsors, etc. 


\section{Orbifold conformal blocks} \label{ConfBlocks}

In this section we review the definition of orbifold conformal blocks
introduced in \cite{FS}.

\subsection{Twisting modules by $\AutDN{N}{}$ }
\label{twistingmodules}

Let $\mc{D}$ be a disk with an $N$--th order automorphism as in the
last section, and let
$\Mt$ be an admissible $\sigma_V$--twisted module over a conformal
vertex algebra $V$. Define a representation $r^{\Mt}$ of the Lie
algebra $\Der^{(o)}_N(\mc{O}) $ on $\Mt$ by the formula
\[
	z^{k+\frac{1}{N}} \partial_{\zf{N}} \rightarrow - N \cdot
	L^{\Mt}_{k}.
\]
It follows from the definition of a twisted module that the operators
$L^{\Mt}_{k}, k > 0$, act locally nilpotently on $\Mt$ and that the
eigenvalues of $L^{\Mt}_{0}$ lie in $\frac{1}{N} \mathbb{Z}$, so that
the operator $N \cdot L^{\Mt}_{0}$ has integer eigenvalues. This
implies that the Lie algebra representation $r^{\Mt}$ may be
exponentiated to a representation $R^{\Mt}$ of the group
$\Aut_{N}(\mc{O})$.  In particular, the subgroup $\Z/N\Z$ of
$\Aut_{N}(\mc{O})$ acts on $\Mt$ by the formula $i \mapsto S_\sigma^i$, where
$S_\sigma$ is the operator defined in \secref{VOA}.

We now twist the module $\Mt$ by the action of $\AutOr{N}$ and define
the vector space
\begin{equation}
	\mathcal{\Mt}(\mathcal{D}) \overset{\on{def}}{=} \AutDN{N}{}
	\underset{\AutOr{N}}\times \Mt. \label{GeoTwistModule}
\end{equation}
Thus, vectors in $\mathcal{\Mt}({\mathcal{D}})$ are pairs $(t,v)$, up
to the equivalence relation $$(\rho(t),v) \sim (t,R^{\Mt}(v)), \qquad
t \in \AutDN{N}{}, v \in \Mt.$$ When $\fd = \fd{x}$, the formal
neighborhood of a point $x$ on an algebraic curve $X$, we will use the
notation $\mc{\Mt}_{x}$.

\subsection{The vector bundle $\V^{G}_{X}$ } \label{VXH}

Let $(X, \bf{p})$ be a smooth $n$--pointed curve over
$\Spec (\mathbb{C})$, and let $\pi: C \mapsto X $ be a smooth $G$--cover of $(X,\bf{p})$. Suppose furthermore that $V$ is a conformal vertex
algebra, and that $G$ acts on $V$ by conformal automorphisms. Let
$\Au_{C}$ be the $\AutO$--torsor over $C$ whose fiber at $y \in C$ is
$\Au(\mc{D}_{y})$the set of formal coordinates at $y$. As explained in \cite{FB},
$\AutO$ acts on $V$, and the action commutes with $G$. Let
\[
\V_{C} = \Au_{C} \underset{\AutO} \times V
\]
The vector bundle $\V_{C}$ carries a $G$--equivariant structure lifting
the action of $G$ on $C$. It is given by
\begin{equation} \label{EqStruct}
	g \cdot (p,(A,z)) \overset{\on{def}}{=} (g(p),(g(A), z \circ
	g^{-1}))
\end{equation}
where $z \circ g^{-1}$ is the coordinate induced at $g(p)$ from
$z$. 
Let $\ol{C} \subset C$ denote the open set on which the $G$--action is
free, and let $\ol{X} \subset X = \pi(\ol{C})$.  Thus, $\ol{C}$ is a
$G$--principal bundle over $\ol{X}$. $\V_{\ol{C}}$ descends to a
vector bundle $\V^{G}_{\ol{X}}$ on $\ol{X}$. More explicitly,
\begin{equation}    \label{wtV}
	\V^{G}_{\ol{X}} = \Au_{\ol{C}} \underset{\AutO \times
	G}{\times} V
\end{equation}
Here, $G$ acts on $\Au_{\ol{C}}$ by $g(p,z)=(g(p),z \circ g^{-1})$,
and this action commutes with the action of $\AutO$. 
The vector bundle $\V^{G}_{\ol{X}}$ possesses a flat connection
$\nabla^{G}$. If $z$ is a local coordinate at $x \in \ol{X}$,
$\nabla^{G}$ is given by the expression $d + L^{V}_{-1} \otimes dz$.

\subsection{Modules along $G$--orbits} \label{OrbitModule}

Let $y \in C$. Then every point $r \in O(y)$ has a cyclic
stabilizer of order $N$, which we denote $G_{r}$. Each $G_{r}$ has a
canonical generator $h_{r}$, which corresponds to the monodromy of a
small loop around $p=\pi(y)$. For a generic point $r$, $N=1$, $G_r = \{ e \}$ and we
set $h_r = e$. Suppose that we are given the following data:

\begin{enumerate}
	\item A collection of admissible $V$--modules $\{ M^{h_{r}}_{r}
	\}_{r \in O(y)}$, one for each point in the orbit, such
	that $M^{h_{r}}_{r}$ is $h_{r}$--twisted.

        \item A collection
	of maps $S_{g,r,g(r)}: M^{h_{r}}_{r} \mapsto
	M^{h_{g(r)}}_{g(r)}$, $g \in G$, $r \in \pi^{-1}(p)$,
	commuting with the action of $\AutOr{N}$ and satisfying
	$$S_{gk,r,gk(r)} = S_{g,k(r),gk(r)} \circ S_{k,r,k(r)},$$
	$$S^{-1}_{g,r,g(r)} = S_{g^{-1},g(r),r},$$ and
	$$S_{g,r,g(r)}^{-1} Y^{M^{h_{g(r)}}_{g(r)}} (g \cdot A,z^{})
	S_{g,r,g(r)} = Y^{M^{h_{r}}_{r}}(A,z).$$
 
	\item If $g \in G_{r}$, then $S_{g,r,r} = S_{g}$, where
	$S_{g}$ is the operator defined in \secref{VOA}.
\end{enumerate}

Given a collection $\{ M^{h_{r}}_{r} \}_{r \in O(y)}$, we can
form the collection $\{ \mc{M}^{h_{r}}_{r}(\fd{r}) \}_{r \in
O(y)}$, where $\mc{M}^{h_{r}}_{r}(\fd{r})$ is the
$\AutOr{N}$--twist of $M^{h_{r}}_{r}$ by the torsor of special
coordinates at $p$. Let 
\[
 \overline{\mc{M}_{O(y)}} = \bigoplus_{r \in O(y)}
 \mc{M}^{h_{r}}_{r}(\fd{r})
\]

This is a representation of $G$, where
$G$ acts as follows. If $A \in M^{h_{r}}_{r}$, $\zf{N}_{r}$ is a
special coordinate at $r$, and $g \in G$, then
\[
g \cdot (A, \zf{N}_{r}) = (S_{g,r,g(r)} \cdot A, \zf{N}_{r} \circ g^{-1}) 
\]
Note that this action is well-defined since the $S$--operators commute
with the action of $\AutOr{N}$. Now, let $\mc{M}_{O(y)} =
(\overline{\mc{M}_{O(y)}})^{G}$, the space of $G$--invariants of
$\overline{\mc{M}_{O(y)}}$. For every $r \in O(y) $ let
\[
\phi_{r} : \mc{M}_{O(y)} \mapsto \mc{M}^{h_{r}}_{r}(\fd{r})
\]
be the isomorphism which is
the composition of the inclusion $\mc{M}_{O(y)}
\to \overline{\mc{M}_{O(y)}}$ and the projection $\overline{\mc{M}_{O(y)}}
\to \mc{M}^{h_{r}}_{r}(\fd{r})$.   

For $v_r \in \mc{M}^{h_{r}}_{r}(\fd{r})$, let
$[v_r]$ denote $ \phi^{-1}(v_{r}) \in \mc{M}_{O(y)}$.  Note that for each
$(A,\zf{N}_{r})_r \in \mc{M}^{h_{r}}_{r}(\fd{r})$, and $g \in G$,
$[(A, \zf{N}_{r})_r] = [(S_{g,r,g(r)} \cdot A, \zf{N}_{r} \circ
g^{-1})_{g(r)}]$ in $\mc{M}_{p}$.

\begin{definition} \label{OrbitModuleDef}
We call $\mc{M}_{O(y)}$ a \emph{$V$--module along $O(y)$}.
\end{definition}

Henceforth, we will suppress the square brackets for
elements of $\mc{M}_{O(y)}$ and refer to $[(A,\zf{N}_{r})_r]$ simply as
$(A,\zf{N}_{r})$.

\subsection{Construction of Modules along $G$--orbits} 
\label{ModuleConstruction}

In this section we wish to describe an induction procedure which yields
a module along $O(y)$ starting with a point $r \in O(y)$ and an
$h_r$--twisted module $M^{h_{r}}_{r}$. Note that when $G_{r}$ is
trivial, this just an ordinary $V$--module $M$.

Thus, suppose we are given $r \in O(y)$, and an
$h_{r}$--twisted module $M^{h_{r}}_{r}$. Observe that the monodromy
generator at the point $g(r)$ is $h_{g(r)} = gh_{r}g^{-1}$, i.e. the
monodromies are conjugate.

\begin{enumerate}
	\item For $g \in G$, define the module
	$M^{gh_{r}g^{-1}}_{g(r)}$ to be $M^{h_{r}}_{r}$ as a vector
	space, with the $V$--module structure given by the vertex
	operator
\begin{equation} \label{NewStructure}
	Y^{M^{gh_{r}g^{-1}}_{g(r)}}(A,z) =
	Y^{M^{h_{r}}_{r}}(g^{-1} \cdot A,z)
\end{equation}
	It is easily checked that this equips
	$M^{gh_{r}g^{-1}}_{g(r)}$ with the structure of a
	$gh_{r}g^{-1}$--twisted module. Furthermore, if $g \in G_{r}$,
	this construction results in an $h_{r}$--twisted
	module isomorphic to $M^{h_{r}}_{r}$.

	\item Recall that $M^{h_{s}}$ is canonically isomorphic to
	$M^{h_{r}}_{r}$ as a vector space by the previous item. Thus,
	if $s \in O(y)$, and $g(s) \ne s$, define
	$S_{g,s,g(s)}$ to be the identity map.

	\item If $g \in G_{s}$, then $g$ is
	conjugate to an element $g' \in G_{r}$. Define $S_{g,s,s} =
	S_{g',r,r}$ also using the canonical identification.

\end{enumerate}

It is easy to check that this construction is well-defined, and
satisfies the requirements of definition \ref{OrbitModuleDef}.

\begin{definition}
We call a module along $O(y)$ obtained via this construction a
\emph{module along $O(y)$ induced from $M^{h_{r}}_{r}$}, and denote it
$\Ind^{O(y)}_{r} (M^{h^{r}}_{r})$.
\end{definition}

\begin{remark} 
If $G_{r}$ is trivial, and $M = V$, then for any $g \in G$, the new
module structure \ref{NewStructure} is isomorphic to the old one, and
so the resulting module $\mc{M}_{p}$ along $O(y)$ is isomorphic
to $\mc{V}^{G}_{p}$, the fiber of the sheaf $\V^{G}$ at $p$.
\end{remark}

\begin{remark}
If $G=G_{r}$, then $r$ is unique, and so any $h_{r}$--twisted module
results in a module along $O(y)$.
\end{remark}

Now, let ${\bf m} = (m_{1}, \cdots, m_{n}) \in G^{n}_{A}$ be a
collection of monodromies. Let $(\pi:C \mapsto X, p_{1}, \cdots,
p_{n}, q_{1}, \cdots, q_{n}) \in \bf{e}^{-1}(\bf{m})$ be a pointed
$G$--cover with the prescribed monodromies. Given a collection $M_{1},
\cdots, M_{n}$ of $V$--modules, such that $M_{i}$ is $m_{i}$--twisted,
the above induction procedure yields a collection 
$\Ind^{O(q_{1})}_{q_{1}}(M_{1}), \cdots, \Ind^{O(q_{n})}_{q_{n}}(M_{n}),$
of modules along $O(q_{i})$. 

\noindent {\bf Note:} We can label $G$--orbits on $C$ by points of $p \in X$. We will use the notation $\mc{M}_{p}$ to denote a module along the $G$--orbit $\pi^{-1}(p)$. 

\subsection{Geometric Vertex Operators}    \label{dual}

Let $p \in X$, and $\mc{M}_{p}$ a $V$--module along $\pi^{-1}(p)$. Let $\mc{D}^{\times}_{p}$ denote the formal punctured disk $\Spec \mc{K}_{p}$.  It
is shown in \cite{FS} that the vertex operator gives rise to a section
$$
\mathcal{Y}^{\mc{M}_{p}, \vee}: \Gamma(\fpd{p},\V^{G}_{\ol{X}} \otimes
\Omega_{\ol{X}}) \to \End {\mc{M}_{p}}.
$$
Moreover, this map factors through the quotient
$$ \label{PointAlgebra}
U(\V^{G}_{p}) \stackrel{\on{def}}{=} \Gamma(\fpd{p},\V^{G}_{\ol{X}}
\otimes \Omega_{\ol{X}})/\on{Im}
\nabla^{G},
$$
which has a natural Lie algebra structure. The corresponding map
\begin{equation} \label{Expandp}
y_{p}:U(\V^{G}_{p}) \to \End {\mc{M}_{p}} 
\end{equation}
is a homomorphism of Lie
algebras. Note that $p$ does not have to lie in $\ol{X}$, but can be
\emph{any} point of $X$.

For each $r \in \pi^{-1}(p)$, composing the above maps with the isomorphism 
$\End \mc{M}_{p} \cong \End \mc{M}^{h_{r}}_{r}(\fd{r})$ induced by $\phi_{r}$ yields a map:
\[
\mathcal{Y}^{\mc{M}_{p}, \vee}_{r}: \Gamma(\fpd{p},\V^{G}_{\ol{X}} \otimes
\Omega_{\ol{X}}) \to \End \mc{M}^{h_{r}}_{r}(\fd{r})
\]
and a Lie algebra homomorphism 
\begin{equation} \label{Expandpr}
y_{p,r}: U(\V^{G}_{p}) \to \End  \mc{M}^{h_{r}}_{r}(\fd{r}).
\end{equation}

\subsection{A sheaf of Lie algebras}    \label{sheaf of la}

Following Section 8.2.5 of \cite{FB}, let us consider the following
complex of sheaves (in the Zariski topology) on $\ol{X}$:
$$
0 \rightarrow \V^{G}_{\ol{X}} \xrightarrow{\nabla}
\V^{G}_{\ol{X}} \otimes \Omega_{\ol{X}} \rightarrow 0
$$
where $\V^{G}_{\ol{X}} \otimes \Omega_{\ol{X}}$ is placed in cohomological
degree $0$ and $\V^{G}_{\ol{X}}$ is placed in cohomological degree
$-1$ (shifted de Rham complex). Let $h(\V^{G}_{\ol{X}})$ denote the sheaf of
the 0th cohomology, assigning to every Zariski open subset $\Sigma
\subset \ol{X}$ the vector space
$$
U_{\Sigma}(\V^{G}_{\ol{X}}) \overset{\on{def}}{=} \Gamma(\Sigma,\V^{G}_{\ol{X}}
\otimes \Omega_{\ol{X}})/\on{Im} \nabla^{G} 
$$
One can show as in Chapter 18 of \cite{FB} that this is a sheaf of
Lie algebras.
 
According to formula \eqref{PointAlgebra}, for any $p \in \Sigma'$,
where $\Sigma' \subset X$ is such that $\Sigma' \cap \ol{X} = \Sigma$,
restriction induces a Lie algebra homomorphism
$\iota_{p}: U_{\Sigma}(\V^{G}_{\ol{X}}) \to U(\V^{G}_{p})$. We denote the image
by $U_{\Sigma}(\V^{G}_{p})$.

\subsection{Conformal Blocks} \label{CB}

Let $(\pi: C \mapsto X, {\bf p})$
be a $G$--cover, and $\mathbb{M} = (\mc{M}_{p_{1}}, \cdots, \mc{M}_{p_{n}})$
a collection of modules along $\pi^{-1}(p_{1}), \cdots, \pi^{-1}(p_{n})$. Let 
\[
\mathbb{F} = \bigotimes \mc{M}_{p_{i}}
\]
Composing the maps \ref{Expandp}  with
the map
\[
U_{X \backslash \bf{p}} (\V^G_{\ol{X}}) \overset{\sum \iota_{p_{i}}}{\rightarrow} \bigoplus^{n}_{i=1}
U(\V^{G}_{p_{i}})
\]
we obtain an action of the Lie algebra $U_{X \backslash {\bf {p}}
}(\V^{G}_{\ol{X}})$ on $\mathbb{F}$.

\begin{definition}    \label{def coinv}
The space of \emph{coinvariants} is the vector space 
$$\mc{H}^{G}_{V}(C,X, {\bf {p}}, \mathbb{M} ) = \mathbb{F} / U_{X
\backslash {\bf {p}}}(\V^{G}_{\ol{X}}) \cdot \mathbb{F} .
$$ The space of \emph{conformal blocks} is its dual: the vector space
of $U_{X \backslash {\bf p }  }(\V^{G}_{\ol{X}})$--invariant
functionals on $\mathbb{F}$ $$C^{G}_{V}(C,X, {\bf p}, \mathbb{M}) = \on{Hom}_{U_{X \backslash {\bf p}
}(\V^{G}_{\ol{X}})}(\mathbb{F}, \mathbb{C}).$$
\end{definition} 

Suppose now that $(\pi: C \mapsto X, {\bf p}, {\bf q})$ is a pointed $G$--cover, that $m_{i}$ denotes the
monodromy at $q_{i}$, and that $M_{1}, \cdots, M_{n}$ is a collection
of $V$--modules such that $M_{i}$ is $m_{i}$--twisted. Then
$\Ind^{O(q_{i})}_{q_{i}}(M_{i})$ is a module along $O(q_{i})$, and we
can apply our definition above with 
\[
\mathbb{M} = ( \on{Ind}^{O(q_{1})}_{q_{1}}(M_{1}), \cdots, \on{Ind}^{O(q_{n})}_{q_{n}}(M_{n}))
\]
In this case, to emphasize the dependence on the points $q_{i}$ in the fiber, we use the notation 
 $\mc{H}^{G}_{V}(C,X, {\bf p, q}, \mathbb{M} )$ and $C^{G}_{V}(C,X,
{\bf p,q}, \mathbb{M})$.

Denote by $\mc{M}_{q_{i}}$ the twist of $M_{i}$ by the torsor of special coordinates at $q_{i}$.   
For each $i$, we have the following commutative diagram of Lie algebras:

\[
\begin{CD}
{ U(\V^{G}_{p_{i}} )} @>{y_{p_{i}}}>>{\on{End} \mc{M}_{p_{i}} }\\
@| @V{\phi_{q_{i}}}VV \\
   U(\V^{G}_{p_{i}}) @>{y_{p_{i},q_{i}}}>>   {\on{End} \mc{M}_{q_{i}} }\\ 
\end{CD}
\]
Letting $\wt{\mathbb{F}} = \bigotimes \mc{M}_{q_{i}}$, the
commutativity of the diagram implies that
\[
\mc{H}^{G}_{V}(C,X, {\bf {p}}, {\bf q}, \mathbb{M} ) \cong
\wt{\mathbb{F}} / U_{X \backslash {\bf {p}}}(\V^{G}_{\ol{X}}) \cdot
\wt{\mathbb{F}} .
\]
and
\[
C^{G}_{V}(C,X, {\bf p},{\bf q}, \mathbb{M}) = \on{Hom}_{U_{X \backslash {\bf p}
}(\V^{G}_{\ol{X}})}(\wt{\mathbb{F}}, \mathbb{C}).
\]


\section{Localization Functors} \label{Localization}

The purpose of this section is to review the general yoga of
Beilinson-Bernstein localization (see \cite{BB1, BB2, BS, FB})
following \cite{FB}.

\begin{definition}
A \emph{Harish-Chandra} pair is a pair $(\g, K)$ where $\g$ is a Lie
algebra and $K$ is an algebraic group, equipped with the following
data: an embedding $\k \subset \g$ of the Lie algebra $\k$ of $K$ into
$\g$, and an action $Ad$ of $K$ on $\g$ compatible with the adjoint
action of $K$ on $\k$ and the action of $\k$ on $\g$. A $(\g,
K)$--\emph{module} is a vector space $V$ carrying compatible actions
of $\g$ and $K$.
\end{definition}

\begin{definition}
Let $Z$ be a variety over $\mathbb{C}$. A $(\g, K)$--action on $Z$ is
an action of $K$ on $Z$, together with a Lie algebroid homomorphism
$\alpha: \g \otimes \mc{O}_{Z} \mapsto T_{Z}$ to the tangent sheaf of
$Z$. The two actions must satisfy the following compatibility
conditions:
\begin{enumerate}
\item The restriction of $\alpha$ to $\k \otimes \mc{O}_{Z}$ is the differential of the $K$--action. 
\item $\alpha(Ad_{k}(a)) = k \alpha(a) k^{-1}$
\end{enumerate}
The action is said to be \emph{transitive} if $\alpha$ is surjective, and \emph{simply transitive} if $\alpha$ is an isomorphism. 
\end{definition}

These definitions extend naturally to the world of pro-algebraic
groups, pro-varieties, and pro-stacks. (see \cite{BB2, BFM}).

Suppose now that $Z \mapsto S$ is a principal $K$--bundle, and that
$Z$ carries a transitive $(\g, K)$--action extending the fibrewise
$K$--action. Let $V$ be a $(\g, K)$--module. The sheaf $V \otimes
\mc{O}_{Z}$ carries an action of the algebroid $\g \otimes
\mc{O}_{Z}$, and it follows from the surjectivity of $\alpha$ that
$V_{stab} = V \otimes \mc{O}_{Z} / ker(\alpha) \cdot (V \otimes
\mc{O}_{Z} )$ is a module for the algebroid $T_{Z}$, and therefore a
$D_{Z}$--module, where the latter denotes the sheaf of differential
operators on $Z$. Moreover, the $K$--equivariance requirement in the
definition of $(\g,K)$--action and $(\g, K)$--module ensures that
$V_{stab}$ is $K$--equivariant, and so descends to a $D_{S}$--module,
which we denote $\Delta(V)$.

We will need a description of the fiber of $\Delta(V)$ at a point $s
\in S$. Let $Z_{s}$ denote the fiber of $Z$ over $s$, and let $\g^{s}$
denote the Lie algebra $ Z_{s} \underset{K}{\times} \g$. The ideals
$ker(\alpha)_{z}, z \in Z_{s}$ give rise to a well-defined Lie ideal
$\g^{s}_{stab} \subset \g^{s}$. Denote $Z_s \underset{K}{\times} V$ by
$\V_{s}$. The action of $\g$ on $V$ induces an action of $\g^{s}$ on
$\V_{s}$.  One can show (see \cite{FB}) that $\Delta(V)_{s} \cong
\V_{s} / \g^{s}_{stab} \cdot \V_{s}$.

\begin{definition}
The functor 
\[
\Delta: \left( (\g,K) - { \on{mod}} \right) \longrightarrow \left( D_{S}- \on{mod} \right)
\]
sending $V$ to $\Delta(V)$ is called the \emph{localization functor}
associated to the $(\g, K)$--action on $Z$.
\end{definition} 

More generally, suppose that $V$ is a module for a Lie algebra $\l$,
which contains $\g$ as a Lie subalgebra and carries a compatible
adjoint $K$--action. $V \otimes \mc{O}_{Z}$ is then a module over the
Lie algebroid $\l \otimes \mc{O}_{Z}$. Suppose also that we are given
a subsheaf of Lie subalgebras $\tilde{\l} \subset \l \otimes
\mc{O}_{Z} $ satisfying the following conditions:

\begin{enumerate}
\item it is preserved by the action of $K$
\item it is preserved by the action of the Lie algebroid $\g \otimes \mc{O}_{Z}$
\item  it contains $ker(\alpha)$
\end{enumerate}

Then, for the same reason as above, the sheaf $V \otimes \mc{O}_{Z} /
\tilde{\l} \cdot ( V \otimes \mc{O}_{Z})$ becomes a $K$--equivariant
$D_{Z}$--module, and so descends to a $D_{S}$--module which we denote
$\widetilde{\Delta}(V)$. Let $s \in S$, and $Z_{s}$ be as
above. Denote by $\l_{s}$ the Lie algebra $Z_{s} \underset{K}{\times}
\l$. The fibers $\tilde{\l}_{z}, z \in Z_{s} $ give rise to a
well-defined subalgebra $\tilde{\l}_{s} \subset \l_{s}$, and the fiber
$\widetilde{\Delta}(V)_{s}$ is isomorphic to $\V_{s} / \tilde{\l}_{s}
\cdot \V_{s}$.

\subsection{Localization of central extensions}  

Suppose as above, that $(\g, K)$ is a Harish-Chandra pair, and that $Z
\mapsto S$ is a $K$--principal bundle with a transitive $(\g,
K)$--action extending the fibrewise $K$--action. Suppose also that
$\g$ has a central extension
\[
0 \mapsto \mathbb{C} \mapsto \ghat \mapsto \g \mapsto 0
\]
which splits over $\k$. Tensoring the above extension by $\mc{O}_{Z}$
we obtain an extension of Lie algebroids
\begin{equation} \label{AlgebroidExtension}
0 \mapsto \mc{O}_{Z} \mapsto \ghat \otimes \mc{O}_{Z} \mapsto \g \otimes \mc{O}_{Z} \mapsto 0
\end{equation}
We will assume henceforth that the extension \ref{AlgebroidExtension}
splits over $ker(\alpha)$, so that we get an embedding of the ideal
$ker(\alpha)$ into $\ghat \otimes \mc{O}_{Z}$. The quotient $\mathcal{T}$ of
$\ghat \otimes \mc{O}_{Z}$ by $ker(\alpha)$ now fits into an extension
of Lie algebroids
\begin{equation}
0 \mapsto \mc{O}_{Z} \mapsto \mathcal{T} \mapsto T_{Z} \mapsto 0
\end{equation}

Let $V$ be a $(\ghat,K)$--module.  The sheaf $V \otimes \mc{O}_{Z}$
carries an action of the Lie algebroid $\ghat \otimes \mc{O}_{Z}$, and
so the quotient $\widehat{V}_{stab} = V \otimes \mc{O}_{Z} /
ker(\alpha) \cdot (V \otimes \mc{O}_{Z})$ carries an action of
$\mathcal{T}$. Let $D'_{Z} = U(\mathcal{T})$, the enveloping algebroid of
$\mathcal{T}$. This is a sheaf of twisted differential operators (a TDO), and
the $\mathcal{T}$--action on $\widehat{V}_{stab}$ extends naturally to an
action of $D'_{Z}$. By the same argument as above, $D'_{Z}$ and
$\widehat{V}_{stab}$ are $K$--equivariant, and so is the action of
$D'_{Z}$ on $\widehat{V}_{stab}$. $D'_{Z}$ therefore descends to a TDO
$D'_{S}$ on $S$, and $\widehat{V}_{stab}$ to a $D'_{S}$--module which
we denote $\Delta(V)$.  Let $\ghat^{s}$ denote the Lie algebra $ Z_{s}
\underset{K}{\times} \ghat$. The ideals $ker(\alpha)_{z}, z \in Z_{s}$
give rise to a well-defined Lie ideal $\ghat^{s}_{stab} \subset
\ghat^{s}$. Denote $Z_s \underset{K}{\times} V$ by $\V_{s}$. The action
of $\ghat$ on $V$ induces an action of $\ghat^{s}$ on $\V_{s}$.  The
fiber $\Delta(V)_{s}$ is isomorphic to $ \V_{s} / \ghat^{s}_{stab}
\cdot \V_{s}$.
Thus, a central extension $\ghat$ of $\g$ gives rise to a localization functor
\[
 \Delta: \left( (\ghat, K)-\on{mod} \right) \longrightarrow \left( D'_{S}-\on{mod} \right) 
\]
sending $V$ to $\Delta(V)$.

More generally, suppose $\ghat$ is a Lie subalgebra of $\widehat{l}$,
and we are given a subsheaf $\widetilde{l} \subset \widehat{l} \otimes
\mc{O}_{Z}$ containing $ker(\alpha)$, preserved by the actions of $K$
and $\widehat{\g} \otimes \mc{O}_{Z}$. Let $V$ be a
$\widehat{l}$--module carrying a compatible $K$--action. The sheaf $V
\otimes \mc{O}_{Z} / \widetilde{l} \cdot (V \otimes \mc{O}_{Z})$ is
then a $K$--equivariant $D'_{Z}$--module, and descends to a
$D'_{S}$--module on $S$, which we denote $\widetilde{\Delta}(V)$. Let
$\widehat{l}_{s} = Z_{s} \underset{K}{\times} \widehat{l}$, and denote
by $\tilde{l}_{s} \subset \widehat{l}_{s}$ the subalgebra arising from
the stabilizers in $Z_{s}$.  The fiber $\widetilde{\Delta}(V)_{s}$ is
isomorphic to $\V_{s} / \tilde{\l}_{s} \cdot \V_{s}$.

\section{$G$--equivariant Virasoro uniformization} \label{Uniformization}

Let $\widehat{\MGg{N}} = \{ (\pi: C \mapsto X, {\bf p}, {\bf q}, {\bf
z}) \}$ where $(\pi: C \mapsto X, {\bf p}, {\bf q})$ is an n- pointed
$G$--cover, and ${\bf z} = (z_1, \cdots, z_n)$, where $z_{i}$ is a
special coordinate at $q_{i}$. This is a projective limit of
Deligne-Mumford stacks. Forgetting the coordinates yields a map
\[
\begin{CD}
{\widehat{\MGg{n}}}\\
@ V{\xi}VV\\
{\MGg{n}}\\
\end{CD}
\]
Similarly, let $\widehat{\Mg{n}} = \{ (X, {\bf p} , {\bf q}, 
{\bf z}') \}$ where $(X, {\bf p}) \in
\Mg{n}$ and ${\bf z}'= (z_{1}^{'}, \cdots, z_{n}^{'})$, where $z^{'}_{i}$ is a formal coordinate at $p_{i}$. Forgetting
the coordinates yields a map
\[
\begin{CD}
{\widehat{\Mg{n}}}\\
@ V{\zeta}VV\\
{\Mg{n}}\\
\end{CD}
\]
Moreover, there exists a map $\eta: \widehat{\MGg{n}} \mapsto
\widehat{\Mg{n}}$ defined by 
$$ \eta: (\pi: C \mapsto X, {\bf p}, {\bf q}, {\bf z}) \mapsto (X, {\bf p},
 {\bf z}')
$$
where $z_{i}'=z_{i}^{N_{i}}$, and $N_{i}$ is the order of $G_{q_{i}}$
making the following diagram commute:
\[
\begin{CD}
{\widehat{\MGg{n}}}@>{\eta}>>{\widehat{\Mg{n}}}\\
@V{\xi}VV@VV{\zeta}V\\
{\MGg{n}}@>>{p}>{\Mg{n}}\\ 
\end{CD}
\]

\noindent For a $\mathbb{C}$-algebra $R$, let 
$$\Aut(C/X, \mc{O}_{R}) = \Aut _{N_{1}}(\mc{O}_R) \times \Aut_{N_{2}} (\mc{O}_{R}) \times \cdots
\times \Aut_{N_{n}}(\mc{O}_R)$$ and 
$$\Der^{(o)}(C/X,\mc{O}) =
\Der^{(o)}_{N_{1}}(\mc{O}) \times \cdots \times \Der^{(o)}_{N_{n}}(\mc{O}).
$$ Similarly, define $\Aut(C/X,\mc{K}_{R})$ and $\Der(C/X, \mc{K})$ in
the obvious way. We have that $\Der(C/X,\mc{O})$
(resp. $\Der(C/X,\mc{K})$) is the Lie algebra of the group scheme
(resp. ind-scheme) representing $R \mapsto \Aut(C/X,\mc{O}_R)$
(resp. $R \mapsto \Aut(C/X,\mc{K}_{R})$).  We see that
$\widehat{\MGg{n}}$ is an $\Aut(C/X)$--torsor over $\MGg{n}$. The following is a generalization of the Virasoro uniformization theorem (see \cite{ADKP}, \cite{BS}, \cite{Kon}, \cite{TUY}) to $G$--covers.

\begin{theorem} \label{GUniformization}
$\widehat{\MGg{N}}$ carries a transitive action of $\Der(C/X, \mc{K})$
extending the action of $\Aut(C/X, \mc{O})$ along the fibers of $\xi$.
\end{theorem}

\begin{proof}
Let $R = \mathbb{C}[\epsilon]/(\epsilon^{2})$. A family of $G$--covers
and special coordinates over $\Spec R$ is the same as a tangent vector
to $\widehat{\MGg{N}}$. Thus to construct an action of $\Der(C/X, \mc{K})$ on
$\widehat{\MGg{N}}$ and prove transitivity, it suffices to do it over $R$.  In fact, we
construct an action of the corresponding group $\Aut(C/X,\mc{K}_R)$. 
Suppose that $(\pi: C \mapsto X,{\bf p }, {\bf q }, {
\bf z })$ is an $R$--point of $\widehat{\MGg{n}}$, and that $\rho \in
\Aut(C/X,\mc{K}_R)$.
Let
$$ \Aut^{G}_{e}(C/X,\mc{K}_R) =
  \Aut^{G}_{e}(\widehat{O(q_{1})}^{*},\mc{K}_R) \times \cdots \times
  \Aut^{G}_{e}(\widehat{O(q_{n})}^{*},\mc{K}_R)
$$
We have an isomorphism 
$$ \psi_{q_{1}}^{*} \times \cdots \times \psi_{q_{n}}^{*}:
\Aut^{G}_{e}(C/X,\mc{K}_R) \mapsto \Aut(C/X,\mc{K}_R ). $$  under
which $\rho$ corresponds to an element $\rho^{G}$ of
$\Aut^{G}_{e}(C/X,\mc{K}_R)$.

At each point $r$ of $O(q_{i})$, there is a set of $N_{i}$ special
coordinates compatible with $z_{i}$ (i.e. obtained from $z_{i}$ by
applying elements of $G$), and we choose one arbitrarily, and
call it $w_{r}$. We now define the new $\rho$--twisted $G$--cover
$$ (\pi: C_{\rho} \mapsto X_{\rho}, {\bf p_{\rho} }, {\bf q_{\rho}
}, {\bf z_{\rho} } )
$$ by the following: as a topological space, $C_{\rho} = C$ and
$X_{\rho} = X$, but the structure sheaf of $C_{\rho}$ is changed as
follows. Let $U \in C$ be Zariski open. If $O(q_{i}) \cap U =
\emptyset, \forall i$, define $\mc{O}_{C_{\rho}}(U) =
\mc{O}_{C}(U)$. If $U$ intersects the orbits $O(q_{i})$ at the points
$r_{1}, \cdots, r_{m}$, define $\mc{O}_{C_{\rho}}(U)$ to be the
subring of $\mc{O}_{C}(U \backslash \{ r_{j} \}_{j=1, \cdots, m})$
consisting of functions $f$ whose expansion $f_{r_{j}}(w_{r_{j}}) \in
R((w_{r_{j}})) $ at $r_{j}$ in the coordinate $w_{r_{j}}$ satisfies
$f_{r_{j}}(\rho^{G,-1}(w_{r_{j}})) \in R[[w_{r_{j}}]]$. Since the gluing is 
$G$--equivariant, $C_{\rho}$ is again a $G$--cover. 

We now prove transitivity. The homomorphisms 
$\mu_{N}: \Aut_{N}(\mc{K}_R) \mapsto \Aut(\mc{K}_R)$
induce a homomorphism
$$ \mu_{C/X} : \Aut(C/X,\mc{K}_R) \mapsto \Aut(X,\mc{K}_R) \overset{def}{=} \Aut(\mc{K}_R) \times
\cdots \times \Aut(\mc{K}_R) \; (n \on{times})
$$
The following diagram commutes:
\[
\begin{CD}
{\widehat{\MGg{n}}}@>{\eta}>>{\widehat{\Mg{n}}}\\
@V{\rho}VV@VV{\mu_{C/X}(\rho)}V\\
{\widehat{\MGg{n}}}@>>{\eta}>{\widehat{\Mg{n}}}\\ 
\end{CD}
\]
The transitivity follows from the fact that the action of
$\Aut(X,\mc{K}_R)$ on $\widehat{\Mg{n}}$ is transitive (by
\cite{ADKP}, \cite{BS}, \cite{Kon}, \cite{TUY}), the fact that $\eta$
is a quasi-finite map, and $\mu_{C/X}$ is surjective.

\end{proof}

Let $O({\bf q}) = \bigcup_{i} O(q_{i})$.
The stabilizer of $(\pi: C \mapsto X, {\bf p }, {\bf q}, {\bf
z })$ is the subgroup $\Aut(C/X,\mc{K})_{out} \in
\Aut(C/X,\mc{K})$ consisting of those elements that preserve
$\mc{O}_{C}(C \backslash O({\bf q}))$. The Lie algebra of
$\Aut(C/X,\mc{K})_{out}$ is the Lie algebra of
$G$--invariant vector fields on $C \backslash O({\bf q})$, denoted 
$\on{Vect}^{G}(C \backslash  O ({\bf q }))$. Note that 
\[
\on{Vect}^{G}(C \backslash  O ({\bf q }) \cong \on{Vect}(X \backslash {\bf p })
\]
where $ \on{Vect}(X \backslash {\bf p }) $ denotes the space of vector
fields on $X \backslash {\bf p }$.
We thus obtain a localization functor
\[
    \Delta:  \left( (\on{Der}(C/X), \Aut(C/X))- \on{mod} \right) \longrightarrow \left( D_{\MGg{n}}- \on{mod} \right)
\]
\[
    \on{M} \mapsto \Delta(\on{M})
\]
The fiber of $\Delta(\on{M})$ at $(\pi: C \mapsto X,{\bf p},{\bf q},
{\bf z})$ is $\mathcal{{M}}/ \on{Vect}^{G}(C \backslash 
O({\bf q}) \}) \cdot {\mc{M}}$, where
$\mathcal{{M}}$ is the $\Aut(C/X)$--twist of $\on{M}$.

We want to construct a localization functor that sends modules along
orbits to $D$--modules whose fibers are spaces of orbifold conformal
blocks. We now proceed with this construction.

\subsection{Construction of the Lie algebroid 
$\mc{U}^{G}(V)_{out}$ }

In this section we construct a Lie algebroid $\mc{U}^{G}(V)_{out}$
over $\widehat{\MGg{n}}$ whose fiber at $(\pi: C \mapsto X, {\bf p},
{\bf q}, {\bf z})$ is the Lie algebra $U_{X \backslash {\bf
p}}(\V^{G}_{\ol{X}})$ of section \ref{sheaf of la}.  Each $G$--cover
in a fixed component $\MGg{n}$ has a fixed genus, which we'll denote
$h$.  This fixed component of $\MGg{n}$, parameterizing a family of
genus $h$ curves with $n$ marked points thus comes with a map to
$\mathfrak{M}_{h,n}$.  Denote by $\mc{X}_{h}$ the universal pointed curve
over $\mathfrak{M}_{h}$. Consider the following fibered product

\[
\begin{CD}
{\widehat{\mc{X}_{h}} =  \widehat{\MGg{n}}} \underset{\mathfrak{M}_{h}}{\times} \mc{X}_{h}@>{}>>{\mc{X}_{h}}\\
@V{}VV@VV{}V\\
{\widehat{\MGg{n}}}@>>{}>{{\mathfrak{M}_{h}}}\\ 
\end{CD}
\]
Thus $\widehat{\mc{X}}_{h}$ can be viewed as the universal pointed
$G$--cover over $\MGg{n}$.  The projection $\widehat{\mc{X}}_{h} \mapsto \widehat{\MGg{n}}$ comes with sections $q_{i}: \widehat{\MGg{n}} \mapsto \widehat{\mc{X}}_{h}$, and applying elements of the group $G$, we obtain sections passing through every point of $O(q_{i})$.

We now wish to construct the coordinate
bundle over $\widehat{\mc{X}}_{h}$. This can be done as follows:
consider the diagram
\[
\begin{CD}
{\mc{S} = \widehat{\MGg{n}} \underset{\mathfrak{M}_{h}}{\times}
\widehat{\mathfrak{M}_{h,1} } } @>{}>>{\widehat{\mathfrak{M}_{h,1}}}\\
@V{}VV@VV{}V\\ {\widehat{\MGg{n}}}@>>{}>{\mathfrak{M}_{h}}\\
\end{CD}
\]
We see that $\mc{S}$ is an $\AutO$--bundle over $\widehat{\mc{X}}_{h}$. Let 
\[
\V_{\widehat{\mc{X}}_{h}} = \mc{S} \underset{\AutO}{\times} V
\]
- this is a sheaf over $\widehat{\mc{X}}_{h}$.
$\V_{\widehat{\mc{X}}_{h}}$ possesses a flat connection $\nabla^{G}$ along the
fibers of the projection $\widehat{\mc{X}}_{h} \mapsto \widehat{\MGg{n}}$, which
along each fiber restricts to the connection $\nabla^{G}$ of section
\ref{VXH}, i.e. we have a complex
\[
0 \mapsto \V_{\widehat{\mc{X}}_{h}} \overset{\nabla^{G}}{\mapsto} \V_{\widehat{\mc{X}}_{h}}
\otimes \Omega_{\widehat{\mc{X}}_{h} / \widehat{\MGg{n}}} \mapsto 0
\]
Restricting this complex to $\widehat{\mc{X}}_{h} \backslash O({\bf q})$, 
taking de Rham cohomology along the fibers, and
then $G$--invariants, yields the Lie algebroid $\mc{U}^{G}(V)_{out}$.

\subsection{Sheaves of Conformal Blocks}

Let ${\bf m} = (m_{1}, \cdots, m_{n})$ be a collection of monodromies,
and $\mathbb{M} = (M_{1}, \cdots, M_{n})$ a collection of $V$--modules
such that $M_{i}$ is $m_{i}$--twisted and admissible. Let $U(M_{i})$
denote the Lie algebra of Fourier coefficients of (twisted) vertex
operators acting on $M_{i}$, and let
\[
U(\mathbb{M}) = U(M_{1}) \times \cdots \times U(M_{n})
\]
The data of the local coordinates in $\widehat{\MGg{n}}$ allows us to embed the Lie algebroid
$\mc{U}^{G}(V)_{out}$ inside $U(\mathbb{M}) \otimes
\mc{O}_{\widehat{\MGg{n}}}$ (by expanding vertex operators in the
local coordinates), and
the Virasoro operators acting on the modules
yield the embedding
\[
\on{Der}(C/X, \mc{K}) \otimes \mc{O}_{\widehat{\MGg{n}}}
\hookrightarrow U(\mathbb{M}) \otimes \mc{O}_{\widehat{\MGg{n}}} 
\]
Denote by $\vect^{G}(C \backslash O({\bf q}))$ the Lie algebroid on
$\widehat{\MGg{n}}$ whose fiber at $(\pi: C \mapsto X, {\bf p}, {\bf
q}, {\bf z})$ is $\on{Vect}^{G}(C \backslash O ({\bf q }))$. Again,
the local coordinates allow us to identify it with a subalgebroid of
$\on{Der}(C/X, \mc{K}) \otimes
\mc{O}_{\widehat{\MGg{n}}}$. Furthermore, under this identification,
\[
\vect^{G}(C \backslash O({\bf q})) \hookrightarrow \mc{U}^{G}(V)_{out}.
\]
The same argument that is used in \cite{FB}, chap. 16, shows that
$\mc{U}^{G}(V)_{out}$ is preserved by the $\Aut(C/X)$-action on
$\widehat{\MGg{n}}$

Let $\on{M} = \bigotimes_{i} M_{i}$. This is a
$U(\mathbb{M})$--module and a $(\g, K)$--module for the Harish-Chandra
pair $(\Der(C/X), \Aut(C/X))$.  We are thus in the setup of section
\ref{Localization}, with $\l =  U(\mathbb{M}) \otimes
\mc{O}_{\widehat{\MGg{n}}} $ and $\wt{\l} = \mc{U}^{G}(V)_{out} $. Suppose first that the
Virasoro central charge $c$ vanishes.  Beilinson-Bernstein
localization then yields a $D$--module $\mc{H}^{G}_{V}(\mathbb{M})$ on $\MGg{n}({\bf m})$ whose
fiber at $(\pi: C \mapsto X, {\bf p}, {\bf q})$ is 
\begin{equation} \label{fiber}
\wt{\mathbb{F}} / U_{X \backslash {\bf {p}}}(\V^{G}_{\ol{X}}) \cdot
\wt{\mathbb{F}}
\end{equation}
where $\wt{\mathbb{F}}$ denotes the $\Aut(C/X)$-twist of $\on{M}$ by
the torsor of special coordinates at the points $q_{i}$. By the
comments at the end of section \ref{CB}, this space is exactly
$\mc{H}^{G}_{V}(C,X, {\bf {p}}, {\bf q}, \mathbb{M} )$.

Suppose now that $c$ is non-zero, and let $\widehat{\Der}(C/X)$ denote
the central extension of $\Der(C/X)$ obtained as the Baer sum of
Virasoro cocycles of the individual factors.  In this case $\on{M}$ is a
Harish-Chandra module for the pair $(\widehat{\Der}(C/X),
\Aut(C/X))$. Moreover, the Lie algebroid extension
\[
0 \mapsto \mc{O}_{\widehat{\MGg{n}}} \mapsto \widehat{\Der}(C/X)
\otimes \mc{O}_{\widehat{\MGg{n}}} \mapsto \Der (C/X) \otimes
\mc{O}_{\widehat{\MGg{n}}} \mapsto 0
\]
splits over the kernel of the anchor map, as the Virasoro cocycle is
trivial on $\on{Vect}^{G}(C \backslash O ({\bf q }))$. Applying the machinery of section
\ref{Localization} with $\widehat{l} = U(\mathbb{M}) \otimes
\mc{O}_{\widehat{\MGg{n}}} $ and $\wt{l} = \mc{U}^{G}(V)_{out}.$, we obtain a module $\mc{H}^{G}_{V}(\mathbb{M})$ for a sheaf of
twisted differential operators on $\MGg{n}({\bf m})$, whose fiber at
$(\pi: C \mapsto X, {\bf p}, {\bf q})$ is $\mc{H}^{G}_{V}(C,X, {\bf
{p}}, {\bf q}, \mathbb{M} )$.

\begin{definition}
Denote by $\nabla^{KZ}$ the connection on $\mc{H}^{G}_{V}(\mathbb{M})$ coming from the $D$--module structure.
\end{definition}

\section{The orbifold KZ connection in the direction of $X$ fixed.}

In this section we give a formula for the orbifold KZ connection along the fibers of the projection
\[
\kappa: \MGg{n}({\bf m}) \mapsto \mathfrak{M}_{g}
\]
sending $(C,X, {\bf p}, {\bf q} )$ to $X$
i.e., this amounts to fixing the base curve $X$. 
We begin by introducing some notation. Fix ${\bf m} = (m_{1}, \cdots, m_{n})$ and $\mathbb{M}=(M_{1}, \cdots, M_{n})$ as before, and let 
\[
\on{M} = M_{1} \otimes \cdots \otimes M_{n}
\]
Let $\Der(\mc{O})$ 
denote the Lie algebra generated by $\{ z^{k} \partial_{z} \}, k \geq 0$, and 
and $\Der_{N}(\mc{O})$ the one generated by $\{ z^{k+\frac{1}{N}} \partial_{\zf{N}} \}, k \geq -1$. 
Under the homomorphism \ref{muN}, these two Lie algebras are isomorphic. 
We have 
\[
\Der^{(o)}(\mc{O}) \subset \Der(\mc{O}) \; \; \; \Der^{(o)}_{N}(\mc{O}) \subset \Der_{N}(\mc{O}) 
\]
Let
\[
\Der(C/X,\mc{O}) =
\Der_{N_{1}}(\mc{O}) \times \cdots \times \Der_{N_{n}}(\mc{O}).
\]
We now wish to consider the two Harish-Chandra pairs $(\Der(\mc{O}), \Aut(\mc{O}))$ and \\ $(\Der_{N}(\mc{O}), \Aut_{N}(\mc{O}))$. 
As explained in \cite{FB}, the action of the pair $(\Der(\mc{O}), \Aut(\mc{O}))$ is simply transitive along the fibers of the projection $\widehat{\Mg{1}} \mapsto \mathfrak{M}_{g}$ (and a similar statement applies in the case of multiple points).  
From the fact that the map
\[
\eta: \widehat{\MGg{n}} \mapsto \widehat{\Mg{n}}
\]
has finite fibers, and that the two Harish-Chandra pairs $(\Der(\mc{O}), \Aut(\mc{O}))$ and \\$(\Der_{N}(\mc{O}), \Aut_{N}(\mc{O}))$ are isogenous, we can deduce that the action of the pair \\
$(\Der(C/X,\mc{O}), \Aut(C/X, \mc{O}))$ along the fibers of the map
\[
\kappa \circ \xi : \widehat{\MGg{n}} \mapsto \mathfrak{M}_{g}
\]
is also simply transitive. Let 
\[
\mathfrak{M}^{G}_{X,n}({\bf m}) = \kappa^{-1}(X) \cap \MGg{n}({\bf m})
\]
and
\[
\widehat{\mathfrak{M}^{G}_{X,n}}({\bf m}) = (\kappa \circ \xi )^{-1}(X) \cap \widehat{\MGg{n}}({\bf m}) 
\]
and let 
\[
\wt{\mc{H}}^{G}_{V}(\mathbb{M},X) = \widehat{\mathfrak{M}^{G}_{X,n}}({\bf m}) \underset{\Aut(C/X)}{\times} \on{M}.
\]
We deduce that along $\mathfrak{M}^{G}_{X,n}({\bf m})$, the connection $\nabla^{KZ}$ on $\mc{H}^{G}_{V}(\mathbb{M})$
lifts to the sheaf $\wt{\mc{H}}^{G}_{V}(\mathbb{M},X)$, and an explicit formula in local coordinates can be obtained. 

Choose coordinates $z_{i}$ around $p_{i}$. Since the map $p$ is quasi-finite, the $z_{i}$ therefore define coordinates on $\kappa^{-1}(X)$. Choose compatible special coordinates $z^{\frac{1}{N_{i}}}_{i}$ at $q_{i}$ which are $N_{i}$--th roots of the $z_{i}$. $z_{i}$ induces coordinates $z_{i} - w_{i}$ at points 
near $p_{i}$, and likewise, the choice of $N_{i}$--th root $z^{\frac{1}{N_{i}}}_{i}$ at $q_{i}$ induces a family of $N_{i}$--th roots of $z_{i}-w_{i}$, which we denote $(z_{i}-w_{i})^{\frac{1}{N_{i}}}$. We can thus trivialize $\wt{\mc{H}}^{G}_{V}(\mathbb{M},X)$ in a neighborhood $W$ of $(C,X,{\bf p}, {\bf q})$:
\[
 \iota: \on{M} \times W \mapsto \wt{\mc{H}}^{G}_{V}(\mathbb{M},X)
\]
\[
(A_{1} \otimes \cdots \otimes A_{n}, w_{1}, \cdots, w_{n}) \mapsto [A_{1}, (z_{1} - w_{1})^{\frac{1}{N_{1}}} ] \otimes \cdots \otimes [A_{n}, (z_{n} - w_{n})^{\frac{1}{N_{n}}} ] 
\]
and similarly, we obtain a trivialization $\iota^{*}$ of  $\wt{\mc{H}}^{G}_{V}(\mathbb{M},X)^{*}$. We have the following theorem:

\begin{theorem}
Along, $\mathfrak{M}^{G}_{X,n}({\bf m})$, the orbifold KZ connection $\nabla^{KZ}$ on $\mc{H}^{G}_{V}(\mathbb{M})$ lifts to $\wt{\mc{H}}^{G}_{V}(\mathbb{M},X)$, and
in the trivialization $\iota$, is given in the local coordinates $z_{i}$ by 
\[
\nabla_{\partial_{z_{i}}} = \partial_{z_{i}} + L^{M_{i}}_{-1}
\]
where $L^{M_{i}}_{-1}$ is the $-1$-st Virasoro operator acting on $M_{i}$.
\end{theorem}
\medskip

\noindent {\bf Note:} The vector field responsible for translations along the curve $X$ is $\partial_{z}$. Observe that under the homomorphism \ref{muN}, we have 
\[
\mu_{N}^{-1} (\partial_{z}) = 1/N z^{-1+1/N} \partial_{\zf{N}}
\]
It follows that $ 1/N z^{-1+1/N} \partial_{\zf{N}}$ is precisely the vector field responsible for moving the ramification points of the map $\pi: C \mapsto X$ along $X$. 
Note that in general, its action on $\widehat{\MGg{n}}$ will change the complex structure of the cover $C$.


\begin{thebibliography}{99}

\bibitem[ACV]{ACV} D. Abramovich, A. Corti, A. Vistoli \emph{Twisted
  bundles and admissible covers}, preprint math.AG/0106211

\bibitem[ADKP]{ADKP} Arbarello, E.; De Concini, C.; Kac, V. G.; Procesi,
C. \emph{Moduli spaces of curves and representation theory.}
Comm. Math. Phys.  117 (1988), no. {\bf 1}, 1-36.

\bibitem[BB1]{BB1} A. Beilinson and J. Bernstein, \emph{Localisation
de $g$-modules}. (French) C. R. Acad. Sci. Paris Sér. I Math.  292
(1981), no. 1, 15--18.

\bibitem[BB2]{BB2} A. Beilinson and J. Bernstein, \emph{A proof of the
Jantzen conjectures}, Advances in Soviet Mathematics {\bf 16}, Part 1,
pp. 1-50, AMS, 1993.

\bibitem[BD]{BD} A. Beilinson, V. Drinfeld, \emph{Chiral Algebras},
Preprint, http://www.math.uchicago.edu/$\sim$benzvi

\bibitem[BFM]{BFM} A. Beilinson, B. Feigin, and B. Mazur,
\emph{Introduction to algebraic field theory on curves}, unpublished
manuscript.

\bibitem[BK]{BK} B. Bakalov, V. G. Kac, \emph{Twisted modules over lattice vertex algebras}, preprint,
math.QA/0402315. 

\bibitem[BS]{BS} A. Beilinson, V. Schechtman, \emph{Determinant
bundles and Virasoro algebras} Comm. Math. Phys. 118 (1988), 651-701.

\bibitem[B]{Borch} R. Borcherds, \emph{Vertex algebras, Kac--Moody
algebras and the monster}, Proc. Natl. Acad. Sci. USA {\bf 83} (1986)
3068-3071.

\bibitem[D]{Dong} C. Dong, \emph{Twisted modules for vertex algebras
associated with even lattices}, J. Algebra \textbf{165} (1994)
91-112.

\bibitem[DLM]{DLM} C. Dong, H. Li and G. Mason, \emph{Twisted
representations of vertex operator algebras}, Math. Ann. {\bf 310}
(1998) 571-600.

\bibitem[FB]{FB} E. Frenkel, D. Ben-Zvi, \emph{Vertex algebras and
algebraic curves}, Mathematical Surveys and Monographs { \bf 88}, AMS,
2001.

\bibitem[FS]{FS} E. Frenkel, M. Szczesny, {\em Twisted modules over
    vertex algebras on algebraic curves}, to appear in Advances in Mathematics.

\bibitem[FFR]{FFR} A. Feingold, I. Frenkel, J. Reis, \emph{Spinor
construction of vertex operator algebras, triality, and $E^{(1)}_8$},
Contemp. Math. {\bf 121}, AMS, 1991.

\bibitem[FLM]{FLM} I. Frenkel, J. Lepowsky, A. Meurman, \emph{Vertex
operator algebras and the monster}, Academic Press, 1988.

\bibitem[FMS]{FMS} D. Friedan, E. Martinec, S. Shenker, {\em Conformal
  invariance, supersymmetry and string theory}, Nuclear Phys. {\bf
  B271} (1986) 93--165.

\bibitem[JKK]{JKK} T. Jarvis, R. Kaufmann, T. Kimura, {\em Pointed
Admissible G-Covers and G-equivariant Cohomological Field Theories},
preprint math.AG/0302316.

\bibitem[K]{K} V. Kac, \emph{Vertex algebras for beginners},
Second Edition, AMS, 1998.

\bibitem[KP]{
KP} V. Kac, D. Peterson, \emph{$112$ constructions of the
basic representation of $E_8$}, in \emph{Anomalies, geometry,
topology}, Agronne, 1985, pp. 276-298, World Scientific, 1985.

\bibitem[Ki]{Ki} A. Kirillov Jr. \emph{On the modular functor associated
with a finite group}, preprint math.QA/0310087.

\bibitem[Kon]{Kon} M. Kontsevich, \emph{ The Virasoro algebra and
Teichmuller spaces}, Funct. Anal. Appl. 21 (1987) no. 2, 156-157.

\bibitem[Li]{LI} H. Li, \emph{Local systems of twisted vertex
operators, vertex operator superalgebras and twisted modules},
Contemp. Math {\bf 193} (1996) 203-236.

\bibitem[Le1]{Lep1} J. Lepowsky, \emph{Calculus of twisted vertex
operators}, Proc. Nat. Acad. Sci. U.S.A. {\bf 82} (1985) 8295--8299.

\bibitem[MP]{MP} J. M. Munoz Porras, F.J. Plaza Martin,
\emph{Automorphism group of $k((z))$: applications to the Bosonic
String}, Comm. Math. Phys, 216 (2001) 609-634.

\bibitem[S]{S} M. Szczesny, \emph{Orbifold chiral correlators from
conformal blocks and KZ equations}, to appear.
	
\bibitem[TUY]{TUY} Tsuchiya, Akihiro; Ueno, Kenji; Yamada, Yasuhiko
\emph{Conformal field theory on universal family of stable curves with gauge
symmetries}.  Integrable systems in quantum field theory and
statistical mechanics, 459--566, Adv. Stud. Pure Math., 19, Academic
Press, Boston, MA, 1989.
\end{thebibliography}
\end{document}